\documentclass{conm-p-l}
\usepackage{amssymb}
\newtheorem{theorem}{Theorem}[section]
\newtheorem{lemma}[theorem]{Lemma}

\newtheorem{question}[theorem]{Question}
\newtheorem{proposition}[theorem]{Proposition}

\theoremstyle{definition}

\theoremstyle{remark}
\newtheorem{remark}[theorem]{Remark}
\numberwithin{equation}{section}

\newcommand {\nn}{{\mathcal{N}}}
\begin{document}

\title{ Measuring sets in infinite groups}
\author{Alexandre V. Borovik}

 \address{Department of Mathematics, UMIST, PO Box
88, Manchester M60 1QD, United Kingdom; {\tt
http://www.ma.umist.ac.uk/avb}}

 \email{borovik@umist.ac.uk}

\thanks{The first author was supported by the Royal Society and
The Leverhulme Trust.}

 \author{Alexei G. Myasnikov}
 \address{Department of Mathematics, The City  College  of New York, New York,
NY 10031, USA; {\tt http://home.att.net/\~~alexeim/index.htm}}
\email{alexeim@att.net}

 \author{Vladimir Shpilrain}
 \address{Department of Mathematics, The City  College  of New York, New York,
NY 10031, USA; {\tt http://zebra.sci.ccny.cuny.edu/web/shpil/}}
\email{shpil@groups.sci.ccny.cuny.edu}
 \date{8 July 2001}


\begin{abstract}
 We are now witnessing a rapid growth of a new part of group theory
which has become known as ``statistical group theory". A typical
result in this area would say something like ``a random element
(or a tuple of elements)
 of a group $G$ has a
property P with probability $p$". The validity of a statement
like that does, of course, heavily depend on how one defines
probability on groups, or, equivalently, how one measures sets in
a group (in particular, in a free group). We hope that
 new  approaches to defining probabilities on groups  
  outlined in this paper create, among other   things,
 an appropriate framework for the study of the
 ``average case" complexity of algorithms on groups.

\end{abstract}

\maketitle

\tableofcontents

\section{Introduction }

 A new part of group theory, often called ``statistical group theory", is
becoming increasingly popular since it connects group theory to other areas
of science, most of all to statistics and to theoretical computer science.

A typical result in this area would say something like ``a random
element (or a tuple of elements) of a group $G$ has a property P
with probability $p$" (see e.g.  \cite{Arzh}, \cite{Fr},
\cite{Olsh}).
 The validity of a statement like that does, of course,
heavily depend on how one defines probability on groups, or,
equivalently, how one measures sets in a group. This is the
problem that we address in the present paper; we feel that it
deserves some discussion at a general methodological level. A much
more technical development of some of the ideas can be found in
\cite{BMR}.

Our approach to statistical group theory is formed by needs of
practical computations with infinite groups. In particular, the
starting point of our study of measures on groups was the desire
to identify measures which can be used in the analysis of the
behaviour of genetic algorithms on infinite groups
\cite{My3,My2,U}. Hence one of our  main requirements for a 
measure is that it reflects the nature of algorithms used for
generating  random or pseudo-random elements of a group. In
Section \ref{sec:conditions} we discuss some other  natural
conditions which would make a measure suitable for practical
computations.

Since  we focus mostly on finitely generated discrete  groups,
almost all  measures defined in the present paper are {\it
atomic}. Recall that a probabilistic  measure on a countable set
is atomic if every subset is measurable. Clearly, the latter
condition is natural in the context of computational group theory
where we are restricted to 
dealing with  finite sets of elements. In Section~\ref{sec:3}, 
we suggest a simple  (perhaps even
naive) general  approach to constructing atomic measures on
countable groups associated with various length (or
``complexity'') functions.

In  Sections \ref{sec:Kolmogorov}, \ref{sec:Kolmogorov2}, we look
at the problem of measuring sets in infinite groups when
pseudorandom elements of a group are generated by a
deterministic process. This naturally leads to Kolmogorov
complexity of words as  a more adequate concept of the ``length
function"  on  a free group. The invariance theorems for
Kolmogorov complexity provide for  a uniform treatment of
different ways to define computational complexity functions on
groups.

Section \ref{se:short-inf} deals with the so-called {\it short
elements bias} which occurs in any generator of random elements of
a group $G$ with a given atomic probability distribution $\mu$.
This effect is based on a simple observation that the measure
$\mu(S)$ of an infinite set $S \subset G$ is essentially  defined
by a few short elements from $S$. There are several principal
approaches to this problem. The most popular one suggests to
consider relative frequencies 
$$\rho_k(S) = \frac{|S\cap B_k|}{|B_k|}$$
(where $B_k$ is the ball of radius $k$ in the Cayley graph of $G$
with respect to a given set of generators) and their behaviour at
``infinity" when $k \rightarrow \infty$.  This leads to the notion
of {\it asymptotic density } of $S$ defined as the following limit
(if it exists):
$$
\rho(S) = \lim_{k \to \infty} \frac{|S\cap B_k|}{|B_k|}.
$$
 Another way
to avoid the short elements bias is to replace a single measure
$\mu$ on $G$ by a parametric family of distributions $\{\mu_L\}$
in which $L$ is  the  average length  of elements  of $G$ relative
to $\mu_L$. In this case the asymptotic behaviour of the set $S$
is described by the limit
$$
\mu_{\infty}(S) = \lim_{L \rightarrow \infty} \mu_L(S).$$
 We discuss asymptotic densities in the last section and  refer to
 \cite{BMR} for a detailed discussion of the second method.

 In Section~\ref{sec:degrees},  we show  how a specific
choice of the probability distribution on a free group allows one
to introduce degrees of polynomial growth ``on average" for
functions on groups. In particular, it is a possible way for
introducing hierarchies of the {\it average case complexity} of
various algorithms for infinite  groups, making meaningful
statements like ``the algorithm works in  cubic
time on average". 

Section~\ref{sec:walks} discusses probability measures generated
by random walks on the Cayley graph of a free group $F_n$. This
class of measures is well suited for analysis of randomized
algorithms on groups and behaves nicely with respect to taking
finite factor groups.

 In  Section  ~\ref{sec:density}, we analyze some limitations of the
 asymptotic density as a tool for measuring sets in infinite groups.
 The most  important one is that it is not sensitive enough, i.e., many
 sets of intuitively different sizes have   the same
 asymptotic density (usually $1$ or $0$). Indeed, in
 \cite{woess}, Woess proved that
 every normal subgroup of infinite index in a free group $F_n$, $n \geqslant 2$,
  has asymptotic density $0$.

Another disappointment is offered by  Theorem~\ref{2.1} which
states that the  set of primitive elements of a free group $F_n$,
$n \geqslant 2$,  has asymptotic   density  $0$. In fact, our proof
provides  lower and upper bounds for the relative
frequencies of the set  of primitive elements of $F_n$, which is a
result of independent interest.  

One way around  this problem is
to consider the so-called {\it growth rate} of the relative
frequency defined by
$$\gamma(S) = \limsup_{k \to \infty}\sqrt[k] {\rho_k(S)}.$$
This is  a  useful characteristic of growth of the set $S$ and
it is  more sensitive than the asymptotic density $\rho(S)$ (see
Section ~\ref{sec:density-2}).

On the other hand, the function  $\gamma(S)$  is not even
additive, which makes it difficult to use as a measuring tool.

  It would be very interesting  to check whether sets having
  the same asymptotic densities (see, for example,
 \cite{Arzh}, \cite{Fr}, \cite{Olsh}) will show  different sizes
 with respect to a more sensitive and adequate measure.

\section{Conditions on a measure}

\label{sec:conditions}

Let  $G$   be a group  (finite or countable infinite).  A
probabilistic measure $\mu$ on $G$ has to satisfy the standard
axioms of a probability space.

Recall that a probability space is a set $X$ together with a
$\sigma$-algebra of subsets $\mathcal{A}$ of $X$ and a
probability measure $\mu : {\mathcal{A}} \longrightarrow R$ which
is  a real-valued function satisfying the following axioms:

\begin{itemize}
\item[{\bf (M1)}]  For any set $S \subseteq {\mathcal{A}}$, ~$\mu(S) \geqslant
0$.

\item[{\bf (M2)}] $\mu(X)=1$.

\item[{\bf (M3)}] If $S_i, i \in I$, is a countable collection of
pairwise  disjoint sets from ${\mathcal{A}}$,   then
$$\mu\left(\bigcup_{i \in I} S_i\right) = \sum_{i \in I} ~\mu(S_i).$$
\end{itemize}

Analysis of behaviour of randomized algorithms on groups requires
estimation of the probability of a given element. A measure $\mu$
on $X$ is  {\em atomic} if it satisfies the following condition

\begin{itemize}
\item[{\bf (M4)}] $X$ is countable and every  subset $S \subseteq X$
is  $\mu$-measurable.
\end{itemize}

In this paper, we shall consider mostly atomic measures on a free
group $F_n$ of rank $n$ with basis $\{x_1, \ldots, x_n\}$. Note
that since every coset with respect to a subgroup $H < F_n$ is
measurable, this defines an induced measure on the factor set
$F_n/H$; this new induced measure is also atomic.

A question that might arise (and which is usually addressed in the
theory of growth of groups) is whether or not a measure should
depend on a particular basis of  $F_n$. Independence of a basis
means invariance of the measure under the action of ${\rm
Aut}\,F_n$. In the case of atomic measures, this implies that
infinitely many singletons $\{g^h\}$, $h\in F_n$, have the same
measure as $\{g\}$ does, which implies $\mu(g) =0$. Hence the only
$({\rm Aut}\,F_n)$-invariant atomic measure on $F_n$ is the singular
measure concentrated at the origin $1$, and  atomic measures
necessarily depend on the choice of a basis in $F_n$.

Finally, we record several  informal  conditions on a measure
which will guide us in evaluating behaviour of various measures:

\begin{itemize}
\item[{\bf (C1)}] a measure  should be natural, i.e., it should meet our
expectations of what the  sizes of various sets in a group are.

\item[{\bf (C2)}] a measure should be sensitive, i.e., sets that
(intuitively) seem to be of different sizes, should have different
measures.

\item[{\bf (C3)}] $\mu(S)$  (or, at least, rather
 tight bounds for  $\mu(S)$) should be  easily computable  for
``natural" sets  $S \subseteq F_n$.

\item[{\bf (C4)}] a measure should admit a natural generator  of random
elements in the group.

\end{itemize}

 In the  rest of the paper we  consider various approaches to constructing 
 measures  on  groups, checking them against our list of conditions
(C1)--(C4).

\section{Atomic probability measures}

\label{sec:3}

 In this section, we  discuss a general method  which allows one to
 define  atomic probabilistic  measures on countable groups.
 Our definition of a measure is based on a notion of ``complexity" of
 elements of  $G$ and a probability distribution  on the set of all
 non-negative integers $\nn$.
 These are  the initial basic objects   which determine the
 intrinsic behaviour of   the corresponding  measure.

Depending on   a problem at hand,  one can use different types of
complexities, for example, it can be  a descriptive complexity, a
computational complexity,  the minimal length of the word
representing  
an  element with respect to a given generating set,
or the length of the normal form of   an element. In general, a
{\em complexity function} (or a {\em complexity}) on a group $G$
is   an arbitrary non-negative integral function $c : G
\longrightarrow {\nn}$ such that for every $n \in \nn$ the preimage
$c^{-1}(n)$ is a finite subset of $G$. Note  that a group $G$
with a complexity function is countable. At the end of the section
we discuss more general complexity functions which allow elements
of $G$ to have infinite complexity and elements with real value
complexities (in this case the group $G$ can be uncountable).

The most important example of a complexity function is the {\it
length function} on a group. Let $G$ be a finitely generated group
with a given finite set of generators $S \subseteq G$. For an
element $g \in G$ by $l_S(g)$ we denote the minimal non-negative
integer $n$ such that $g = y_1 \ldots y_n$ for some $y_i \in S
\cup S^{-1}$. If $H$ is a subgroup of $G$ then the restriction of
$l_S$  on $H$ gives rise to a new complexity function on $H$,
which might not be a length function on $H$.

 In Section \ref{sec:Kolmogorov} we consider another important type of
  complexity functions on $G$, so-called Kolmogorov complexity
  functions.

For a given complexity function $c$ on $G$ and a discrete
probability distribution on non-negative integers $\nn$,
given by a density function $d : \nn \longrightarrow R$,    we are
going to construct an atomic measure $\mu_{c,d}$  on $G$.

Recall that, for an  atomic measure on a countable set $X$, the
value of $\mu$ on an arbitrary subset $S \subset X$ is defined
uniquely by the formula
\begin{equation}
\label{eq:mu}
 \mu(S) = \sum_{w \in S} \mu(w).
\end{equation}

This shows that an atomic measure $\mu$ is completely determined
by its values on singletons $\{w\}, \ w \in X.$ In other words, to
define $\mu$ it suffices to define a function $p:X \longrightarrow
R$ (and put $\mu(w) = p(w)$), which is called a {\it probability
mass function} or a {\em density function} on $X$, such that:
$$ p(w) \geqslant 0 \quad \hbox{ for all }\quad  x \in X,$$
$$\sum_{w \in X}p(w) = 1.$$
Now we put forward one more condition on the measure $\mu$ which
ties $\mu$ to a complexity function $c: G \rightarrow 
\mathcal{N}$:
\begin{itemize}
\item[{\bf (C5)}] for any $u, v \in G$,  $c(u) = c(v)$ implies $\mu(u) =
\mu(v)$.
\end{itemize}
We call a measure $\mu$ on $G$ {\it $c$-invariant} if $\mu$
satisfies  (C5), i.e., elements of the same $c$-complexity
have the same measure.

The following result describes $c$-invariant measures on $G$. For
$k \in \mathcal{N}$ define the $k$-sphere $C_k$ with respect to a
complexity $c$ as follows:
$$
C_k  = \{w \in G \mid  c(w) = k\}.
$$
Similarly, by $B_n$ we denote the disc or the ball of radius $n$
with respect to $c$:
 $$B_n = \{w \in G \mid c(w) \leqslant n\}.$$

\begin{lemma}
Let $\mu$ be an atomic measure on $G$ and $c$  a complexity
function on $G$. Then:
  \begin{enumerate}
\item [(1)] If $\mu$ is $c$-invariant, then the function $d_{\mu}:
\mathcal{N} \longrightarrow R$ defined by
$$d_{\mu}: k \longrightarrow \mu(C_k)$$
is a probability measure on $\mathcal{N}$;
\item [(2)] if $d:\mathcal{N} \longrightarrow R$ is a probability
measure on $\mathcal{N}$, then the function $p_{c,d} : G
\longrightarrow \mathcal{R}$  defined by
$$p_{c,d}(w) = \frac{d(c(w))}{\mid C_n \mid}$$
is a probability function which gives rise to an  atomic
$c$-invariant measure on $G$.
\end{enumerate}
\end{lemma}

\medskip
\noindent
{\it Proof\/} is obvious. \hfill $\square$

\begin{remark}
If $c(w)$ is the length of an element $w \in G$ with respect to a
given finite set of generators of $G$, then   $c$-invariant
measures play an important role in asymptotic group theory. We
will call such measures {\it homogeneous}.
\end{remark}

Since $d(k) \longrightarrow 0$ as $k \longrightarrow \infty$, the
following elementary but fundamental property of $\mu_{c,d}$
holds.

\smallskip
\begin{quote}
\noindent  {\it  More complex (with respect to a given complexity
$c$) elements of a set $S \subset X$ contribute less to
$\mu_{c,d}(S).$}
\end{quote}
\smallskip

This  discussion shows that   in defining the probability
measures $\mu_{c,d}$ on $G$, everything  boils down to the problem
of choosing a complexity function $c:G \longrightarrow \nn$ and a
moderating distribution  $d:\nn \longrightarrow \mathbb{R}$. This
 choice might  depend on a particular problem one would like to
address.  To that end, it seems reasonable to:

\begin{itemize}
\item[{\bf (a)}] use  complexity functions on $G$ which reflect
intrinsic features of the problem at hand;

\item[{\bf (b)}]  use   probability distributions on $\nn$ which reflect
the nature of the process which generates (pseudo)random elements
of the group, or which allow to study the statistical properties
of the group in a meaningful and computationally feasible way.
\end{itemize}

We discuss possible complexity functions on groups in the next
section.   Here are some well-established parametric families of
density functions on $\nn$ that we have  considered in this framework:
\smallskip

\noindent (i) The {\it Poisson density}:
$$
d_{\lambda}(k) = e^{-\lambda} \frac{\lambda^k}{k!}
$$

\smallskip
\noindent (ii) The {\it exponential density}:
$$
d_{\lambda}(k) = (1-e^{-\lambda}) e^{-\lambda k}
$$

\smallskip
\noindent (iii) The {\it standard normal\/ {\rm (}Gauss{\rm )} density}:
$$
d_{a, \sigma}(k) = b \cdot  e^{-\frac{(k-a)^2}{\sigma}}
$$

\smallskip
\noindent  (iv) The {\it Cauchy density}:
$$
d_{\lambda}(k) = b \cdot \frac{1}{(k-\lambda)^2 +1}
$$

\smallskip
\noindent (v) The {\it Dirac density}:

 \[ d_m(k) = \left\{\begin{array}{ll}
                1 & \mbox{if $k = m$}\\
                0 & \mbox{if $k \neq m$}
                \end{array}
             \right. \]

The following density function depends on a given complexity
function $c$ on $G$:
\smallskip \noindent (vi) The {\it finite disc uniform density}:

 \[ d_m(k) = \left\{\begin{array}{ll}
                \frac{|C_k|}{|B_m|} & \mbox{if $k \leq m$}\\
                0 & \mbox{otherwise}
                \end{array}
             \right. \]

For example, the Cauchy density functions  are very convenient for
defining degrees of polynomial growth ``on average'', while the
exponential density functions are suitable for defining degrees of
exponential growth ``on average" (see Section~\ref{sec:degrees}).

On the other hand, different distributions arise from different
random generators of elements in a group.  For example, Dirac
densities correspond to the uniform random generators on the
spheres $C_n$; finite disc densities arise from uniform random
generators on the discs $B_n$.  Measures on the free group $F_n$
associated with the exponential moderating distribution on $\nn$
have especially nice properties and are studied in some detail in
\cite{BMR}.

\begin{remark} Sometimes we allow for some elements in the set $X$ to
have infinite complexity $\infty$. In this case we consider
complexity functions of the type $c:X \longrightarrow
N\cup\{\infty\}$ such that $c^{-1}(n)$ is a finite subset of $X$
for every integer $n$ (so that we may have infinitely many
elements in $X$ with complexity $\infty$). Following our
requirement that more complex elements contribute less to the
measure, we define the probability density function $p$  on $X$
  for elements with finite complexity, and for  elements
$x \in X$ of infinite complexity  we just put $p(x) = 0$.
\end{remark}

\section{Kolmogorov complexity functions}

\label{sec:Kolmogorov}

\subsection{Kolmogorov complexity functions on free groups}

In this section we discuss complexity functions on free groups. We
are going to represent elements of free groups as reduced words in
a given basis,  so the starting point of this discussion is
complexity of {\it words}.

Let $A$ be a finite alphabet and  $A^*$  the  set of all finite
words in $A$. The length $|w|$ of  a word  $w = x_1\ldots x_n,$  $x_i
\in A$, is equal to $n$.  The length function
$l:w \longrightarrow |w|$ maps $A^*$ into $\nn$. Observe that
$l^{-1}(n) = |A|^n$. It follows that

\begin{itemize}
\item The length function $l:A^* \longrightarrow \nn$ is a complexity
function on $A^*$.
\end{itemize}
This is one of the basic complexity functions that we will
consider in this paper.  The problem however is that some very long words in
$A^*$ do not look very complex. For example, the word
$a^{10^{100}}$, where $a \in A,$ has length $10^{100}$, but this word is
very easy to describe. This leads to the notion of {\em
Kolmogorov} complexity, or {\em descriptive} complexity, or
sometimes it is called {\em algorithmic} complexity. We refer to
\cite{LV}, \cite{W}  for detailed treatment of Kolmogorov
complexity.

Intuitively,  Kolmogorov complexity of a word $w \in A^*$ is
the minimum possible size (length)  of a   description of $w$ with respect to
a given formal general  procedure. For example, one may think of
the descriptive complexity of $w$ as of the the minimum possible size
  of a  program, in a given programming language, which produces $w$ after
 finitely many steps. In this case, the word $a^{10^{100}}$ above
will have complexity much less than $10^{100}$.

 Now we give a
formal definition of Kolmogorov complexity.

It is sufficient to consider programs of a very particular type, say,
Turing machines.  Denote by $B = \{0,1\}$ the standard binary
alphabet and by $B^*$ the set of all finite binary strings. Every
Turing machine $M$ determines  a  partial (perhaps empty)  recursive
function $f_M: B^* \longrightarrow A^*$. The function $f_M$ is
defined on $x \in B^*$ if and only if the machine $M$ starts on
the tape with a word $x$, halts in finitely many steps, and a
word $w \in A^*$ is written on the tape. In this event, $f_M(x) =
w$.
 Denote by $s_M(x)$ and $t_M(x)$, respectively,  the tape space
and the number of steps  needed for $M$ to write $w$ and  halt.

The Kolmogorov complexity of a word $w \in A^*$ with respect to a
given machine $M$  is  defined as follows:
\[ K_M(w) = \left\{ \begin{array}{ll}
\min\{|x| \mid x \in B^* \ \& \ f_M(x) = w\}  & \mbox{if such $x$
exists}\\ \infty  & \mbox{otherwise}
\end{array}
\right. \]

Similarly, one can define {\em time bounded} and {\em space
bounded } Kolmogorov complexities of a word $w \in A^*$.   Namely,
for a Turing machine $M$  and a  function $$\beta: \nn
\longrightarrow \nn$$ (a {\em bound}) define the space bounded and
time bounded Kolmogorov complexities $KS_M(x,\beta)$ and
$KT_M(x,\beta)$ (with respect to the bound $\beta$) as follows:
 \[ KS_M(w, \beta) = \left\{ \begin{array}{ll}
\min\{|x| \mid x \in B^* \ \& \ f_i(x) = w  \  \& \ s_M(x)
\leqslant \beta(|x|)\}  & \mbox{if such $x$ exists}\\ \infty  &
\mbox{otherwise}
\end{array}
\right. \]

 \[ KT_M(w, \beta) = \left\{ \begin{array}{ll}
\min\{|x| \mid x \in B^* \ \& \ f_i(x) = w  \  \& \ t_M(x)
\leqslant \beta(|x|)\}  & \mbox{if such $x$ exists}\\ \infty  &
\mbox{otherwise}
\end{array}
\right. \]

These definitions depend on  a given Turing machine $M$. It turns
out that  a {\em universal} Turing machine $M$ provides an optimal
notion of Kolmogorov complexity. Namely, the following invariance
theorem holds.

\medskip \noindent {\bf Invariance Theorem I (\cite{Sol}
\cite{Kol}, \cite{Cha}).} {\it There exists a Turing machine $U$
such that for any Turing machine $M$ and  for any word  $w \in
A^*$ the following inequality holds: $$ K_U(w) \leqslant K_M(w) +
C_M,$$  where $C_M$ is a constant which does not depend on $w$.}

\medskip

Similar results hold for space bounded and time bounded
Kolmogorov complexities.

\medskip \noindent {\bf Invariance Theorem II \cite{Har}.}
{\it There exists a Turing machine $U$ such that for any Turing
machine $M$, for any bound $\beta: \nn \longrightarrow \nn$,  and  for
any word  $w \in A^*$ the following inequalities  hold:
 $$ KS_U(w, C_M\beta) \leqslant K_M(w, \beta) + C_M,$$
 $$ KT_U(w, C_M\beta \log(\beta) \leqslant K_M(w, \beta) + C_M,$$
  where  $C_M$ is a constant which does not depend on $w$.}

\medskip

 The theorems above allow one to consider $K_U(w)$ as  the  Kolmogorov complexity
 $K(w)$ of a word $w \in A^*$, were $U$ is  an arbitrary  fixed optimal
 machine.

 We will use these different types of Kolmogorov complexities in constructing
  measures on a free group. But first,  two  remarks are in
  order (a bad news and a good news).

\begin{remark} The function $ w \longrightarrow K(w)$
 is not recursive.
 Thus, we cannot effectively compute the Kolmogorov complexity of  a
given word. However, it turns out that we can estimate $K(w)$ by the length
of $w.$ This shows that the length $|w|$ as a complexity of $w$ is
back in the game.
\end{remark}

\begin{remark} There exists a constant $C$ such that for any word $w \in A^*$
$$ K(w) \leqslant |w|\log_2|A| + C.$$
Actually, there are  better estimates which will be discussed
elsewhere.

\end{remark}

Now we define Kolmogorov complexity of an element of a free group
$F = F(X)$ with a basis $X$. Let $X^{-1} = \{x^{-1} \mid x \in X\}$
and  $ A = X \cup X^{-1}.$ For an element $f \in F$ denote by
$w_f$ the unique
 reduced word in the alphabet $A$ which represents $f$.

  Now  the
 Kolmogorov complexity $K(f,X)$ of $f$, with respect to the basis $X$, is defined as
 $$K(f) = K(w_f).$$

 It is easy to see that for different bases $X$ and $Y$ of $F$, the
corresponding
 Kolmogorov complexities  are equivalent up to some additive  constants, i.e.,
 there are positive integers $C_1$ and $C_2$ such that for any $f \in F$,
 $$
K(f,Y) - C_1 \leqslant  K(f,X) \leqslant K(f,Y) + C_2. $$

This allows us to fix an arbitrary basis of the free  group $F$
and consider all the Kolmogorov complexities with respect to this
particular  basis.

\section{Kolmogorov complexity functions on finitely generated groups}

\label{sec:Kolmogorov2}

Let $G$ be a group generated by a finite set $X$. There are several
 ways of introducing Kolmogorov complexity on $G$.
\smallskip

{\bf Method I.} Let $F(X)$ be a free group on $X$ and let $\eta :
F(X) \longrightarrow G$ be the canonical epimorphism.  For an
element $g \in G$ define the Kolmogorov complexity $K(g,X)$ of
$g$ with respect to $X$ as follows
 $$ K(g,X) = \min\{K(w,X) \mid \eta(w) = g \},$$
 i.e., the Kolmogorov complexity of $g \in G$ is the minimum of
 Kolmogorov complexities of representatives of $g$ with respect to
 $X$.

As we have   already seen, $K(g,X)$ can be estimated from above  by
the geodesic length of $g$ in the Cayley graph of $G$ with respect
to the set $X$ of generators.   It is extremely difficult to deal
with this type of complexity if  the set of such geodesics is itself
  very complex. In this event, it might be useful to
consider some special normal forms of   elements. For example,
this is   the case for  braid groups  with respect to the set
of  Artin generators.
\smallskip

{\bf Method II.} Suppose the group $G$ has a set of normal forms
$V$, i.e.,  there exists a subset $V \subset F$ such that $\eta
\mid_V$ is one-to-one.  Then one can define Kolmogorov  complexity
of $g \in G$ with respect to $V$ as
 $$
 K(g,V) = K(\eta^{-1} |_V(g)).
 $$

 The following method provides an {\em average} Kolmogorov  complexity
 with respect to a given measure on $F$.
\smallskip

{\bf Method III.} Let $\eta :F \longrightarrow G$ be as above.
Let $d: \nn \longrightarrow R$ be a probability  on $\nn$
 and $c: F \longrightarrow \nn$  a Kolmogorov complexity function on
$F$. As we have discussed above, there exists a probability measure $\mu =
\mu_{c,d}$ on the group $F$. Define an {\em average} Kolmogorov
complexity of $g \in G$ with respect to $d$, $c$,  and $X$ as
follows:
 $$ KA(g,d,c,X) =
 \frac{\sum_{w \in \eta^{-1}(g)} c(w)\mu_{c,d}(w)}{\sum_{w \in \eta^{-1}(g)}\mu_{c,d}(w)}.
 $$

\section{Short elements bias and behaviour at infinity}
\label{se:short-inf}

In this section we discuss how to deal with short elements bias in
an infinite group.

Let $G$ be a finitely generated group with a finite set $X$ of
generators. For simplicity we will discuss only measures
corresponding to the length function $l_X$ on $G$, but  similar
arguments can be applied for other complexity functions as well. Let
$$C_n = \{g \in G \mid l_X(g) = n\}, \ \ B_n = \{g \in G \mid l_X(g)
\leqslant n\}$$

Let $\mu$ be an atomic function on $G$.
 Since $\sum  \mu(C_k) =1$, we have 
 $\mu(C_k) \longrightarrow 0$ as $k  \longrightarrow \infty$.
 Therefore elements  of bigger length in a
set $S \subseteq G$ contribute  less to $\mu(S)$, i.e., we witness
the {\it short elements bias}.  On one hand, this meets our
intuitive expectations because, technically, random elements of a
very big length are inaccessible to us, e.g., no computer can
generate for us a random element of length $> 10^{100}$. On the
other hand, it may happen that only a few short elements
essentially define the measure of an infinite set.

There are several approaches to deal with the  short elements bias.

{\bf Method I.}

 In practical computations with groups, when we 
wish to evaluate performance   of a given algorithm $\mathcal{A}$,
a typical solution to this problem is the following. Choose a
sufficiently big random positive integer $n$ (or several of them),
generate (pseudo) randomly and uniformly enough elements of length
$n$, and run your algorithm on the produced inputs. The choice of
$n$ usually  depends on the computer resources and the hardness of
the algorithm. Mathematically this can be  modelled by the
probability distribution $\mu_{l_X,d_n}$ with Dirac density $d_n$.
The only problem is the choice of $n$. Theoretically, to avoid the
bias toward short elements, one has to take the limit when $n
\rightarrow \infty$. More precisely, let $R$ be a subset of $G$.
Denote by $\rho_n(R)$ the measure  of $R$ with respect to Dirac
density concentrating at $n$:
$$
\rho^{(s)}_n(R) = \frac{|R \cap C_n|}{|C_n|}.$$ Then the
asymptotic behaviour of the set $R$ (with respect to Dirac
densities) can be characterized by the following limit (if it
exists):
$$
\rho^{(s)}(R) = \lim_{n \rightarrow \infty} \rho^{(s)}_n(R)$$ This
limit is called the {\it shperical asymptotic density} of the set
$R$.

Similarly, if we measure $R$ with respect to the disc uniform
density function (with respect to the complexity function $l_X$):
$$
\rho^{(d)}_n(R) = \frac{|R \cap B_n|}{|B_n|},$$ and then take the
limit (if it exists)
 $$
\rho^{(d)}(R) = \lim_{n \rightarrow \infty} \rho^{(d)}_n(R)$$ then
we get the {\it disc asymptotic density} of the set $R$.

In Section \ref{sec:density} we discuss asymptotic densities in
detail, here it is worthwhile to mention only that these
characteristics are not sensitive enough to distinguish various
subsets of groups. For example, all subgroups of infinite index
have the same asymptotic density (disc or spherical) equal to $0$.

{\bf Method II.}

 Let $\mu$ be a probability distribution on $G$. Then the {\it
 mean length} $L_{\mu, X}$ of elements of $G$ with respect to $\mu$ and the set
$X$  of generators is the expected value of the  length function $l_X$:
 $$L_{\mu,X} = \sum_{g \in G} l_X(g)\mu(g).$$

For example, if $G = F(X)$ and a measure $\mu_{\lambda} =
\mu_{l_X,d_{\lambda}}$ is given on $F(X)$ by the exponential
density:
$$
d_{\lambda}(k) = (1-e^{-\lambda}) e^{-\lambda k}
$$
then the mean length $L_{\lambda}$ of words in $F$ distributed
according to $\mu_{\lambda}$ is equal to
 $$L_{\lambda} = \sum_{w \in F} |w|\mu(w)
= (1-e^{-\lambda})\sum_{k= 1}^{\infty}k(e^{-\lambda})^{k} =$$
$$= \frac{(1-e^{-\lambda})}{e^{\lambda}} \sum_{k=
1}^{\infty}k(e^{-\lambda})^{k-1}=  \frac{1}{(e^\lambda
-1)e^\lambda}.$$

Hence $L_{\lambda} \rightarrow \infty$ when $\lambda \rightarrow
0$. Thus,   we have

\begin{quote}
{\it a family of probabilistic distributions $\{\mu
_\lambda\}_{\lambda}$ with parameter $\lambda$ such that the mean
length $L_\lambda$ tends to $\infty$ when the parameter $\lambda$
approaches $0$.}
\end{quote}

Similar results hold for all other distributions introduced in
Section \ref{sec:3} (after renormalizing the parameters).

Now one can measure the behaviour of $R$ at infinity by the
following limit (if it exists)
 $$
 \mu_{\infty}(R) = \lim_{\lambda \rightarrow 0}\mu_\lambda(R).$$

A study of  limits of this type was initiated in \cite{BMR}.

Note that if $\{\mu_k\}_{k}$ is the Dirac parametric family of
distributions, then
$$ \mu_{\infty}(R) = \rho^{(s)}(R),$$
 and if $\{\mu_k\}_{k}$ is the parametric family of disc uniform
distributions, then

$$ \mu_{\infty}(R) = \rho^{(d)}(R).$$

\section{Degrees of polynomial growth ``on average"}

\label{sec:degrees}

Let $c: F_n \longrightarrow {\mathbb{R}}$ be a complexity function
and $\mu$ the measure moderated by the Cauchy distribution
$$d(k) = \frac{6}{\pi^2}\cdot \frac{1}{k^2}.$$
The advantage of this distribution is that it allows to measure
degrees of polynomial growth of functions on $F_n$ ``on average''
in the following sense.

Let  $f: F_n \longrightarrow {\mathbb{R}}$ be a non-negative real
valued function. We say that $f$ has polynomial growth of degree
$\alpha \geqslant 0$ if $\alpha$ is the greatest lower bound of
the set of real positive numbers $\beta$ such that the mean value
of $\frac{f(w)}{c(w)^{\beta-1}}$ is finite, that is,
\begin{equation}
\alpha = \inf \left\{\, \beta \;\left|\; \sum_{w\in
F_n}\frac{f(w)}{c(w)^{\beta-1}} \mu_{c,d}(w) \;\; \hbox{
converges}\,\right.\right\}. \label{eq:degree}
\end{equation}
It follows  immediately  from the construction of our measure that
the function $c(w)^m$ has   growth of degree $m$.  (Just
  recall  that the series $\sum
\frac{1}{n^{1+\epsilon}}$ converges for all $\epsilon > 0$, while
the harmonic series $\sum \frac{1}{n}$ diverges.)

In particular, if $f(x)$ is the running time of some algorithm
with input $x$, this definition allows us to make meaningful
statements like ``the algorithm works in  cubic
time on average''.

Now we are going to generalize this situation and try to find sufficient
conditions  to define the degrees of growth of
functions on an arbitrary infinite factor group of $F_n$.

Let $G = F_n/R$ be an infinite factor group of the free group
$F_n$ and $\eta: F_n \longrightarrow G$ the natural homomorphism.
We shall list some natural conditions for the complexity function
$c$ on $F_n$ and the moderating probability distribution $d$
which allow us to define degrees of the average growth of functions
on $G$.

\begin{itemize}
\item[{\bf (D1)}] ({\em Existence of a complexity function on $G$.}) The mean
complexity
$$
\bar c(\bar g) = \frac{\sum_{w \in \bar g}
c(w){\mu_{c,d}}(w)}{\sum_{w \in \bar g} {\mu_{c,d}}(w)}
$$
is finite for every $\bar g \in G$.
\end{itemize}

\begin{itemize}
\item[{\bf (D2)}] ({\em Existence of degrees}.) There exists a positive
number $\tau$ such that the series
$$
\sum_{\bar g \in G} \bar c(\bar g)^a \bar\mu(\bar g)
$$
converges for $a < \tau$ and diverges for $a > \tau$; $\bar \mu$
here is the measure on $G$ induced by the measure $\mu_{c,d}$ on
$F_n$.
\end{itemize}
If these conditions are satisfied, we can define the degree of
growth of an arbitrary nonnegative  real valued function $f(\bar
g)$ on $G$ as
\begin{equation}
(\hbox{degree of growth of } f) = \inf \left\{\,\beta
\;\left|\;\;\sum_{\bar g \in G}\frac{f(\bar g)}{\bar c(\bar
g)^{\beta-\tau}} \bar\mu(\bar g)\;\;
\hbox{converges}\,\right.\right\}. \label{eq:degree3}
\end{equation}
Now we at least have the (easy to check) property that the degree
of growth of the function $\bar c(\bar g)^m$ is $m$.

Note that not every moderating probability distribution $d$ on
$\mathbb{N} ~\cup \{0\}$ is suitable for defining  degrees of polynomial
growth. For example, if we take in the definition  of the
degrees of growth on the free group $F_n$, the exponential
distribution
$$
d(k) = (1 - e^{-\lambda})\cdot e^{-\lambda k},
$$
we  note that the series
$$
\sum_{x \in F_n} c(x)^a \mu_{c,d}(x)
$$
converges for all $a$.

\begin{question}
Can one find a moderating probability distribution\/ $d$ such that
the conditions\/ {\rm D1} and\/ {\rm D2} are satisfied for every
factor group\/ $G = F_n/R$ of infinite index and the measure
$\bar\mu$ on\/ $G$ induced from $\mu_{c,d}$?
\end{question}

The  Cauchy distribution $d(k) = \frac{6}{\pi^2} \cdot
\frac{1}{k^2}$ is still on the list of candidates for the
affirmative answer.

We need to warn  that the degree of growth of a function $f$ on
$G$  is not necessarily equal to the degree of growth of its lift
$f \circ \eta: F_n \longrightarrow \mathbb{R}$. It might happen
that, for certain problems, it is much more convenient to work in
  a free group than in its factor group $G$ and evaluate the
degrees of growth of the lifted function $f \circ \eta$ instead
of those of $f$. Therefore general results concerning relations between
degrees of growth on $F_n$ and $G = F_n/R$ would be rather interesting.

\section{Measures generated by random walks}

\label{sec:walks}

\subsection{Random walks} An interesting and, in some cases,
easier to analyze class of measures is  related to random walks
on   Cayley graphs of groups. From   an algorithmic  point of
view, the most natural way to produce a random element in a
finitely generated group is to first make  a random word on the
generators, and then apply relations.
  This is, in disguise,  a random walk on the
Cayley graph of a given group. This becomes especially relevant when
hardware random numbers generators are used  to produce random words.

Let us look at this basic procedure in more detail.

Let $F_n$ be a free group on  free generators $x_1, \ldots,
x_n$. We can associate with it a free monoid $M_n$ with the
generators $x_1, X_1, \ldots, x_n, X_n$ and the natural
homomorphism $\pi: M_n \longrightarrow F_n$ which sends $x_i$ to
$x_i$ and $X_i$ to $x_i^{-1}$. A random walk of length $l$ on the
Cayley tree of $F_k$ which starts at the identity $1$ is
naturally described by a word of length $l$ in $M_n$.

In essence, we take $M_n$, not $F_n$, as the ambient algebraic
structure, and introduce measures and complexity functions on
$M_n$ rather than on $F_n$. There is a compelling evidence that, in
at least some problems, it might be convenient to work in the
ambient free monoid $M_n$. For example, it appears that
physicists prefer to use the random walk approach in their study
of statistics of braids, knots and heaps, the latter being
closely related to  {\em locally free groups}
$$
{\mathcal{LF}}_n = \langle f_1,\ldots, f_{n} \mid [f_i,f_j] = 1
\hbox{ for } |i-j| > 1\rangle
$$
(see e.g. \cite{CN,DN1,DN2,VNB}). In any case, the physical
process of the accumulation of soot (in the 2-dimensional case) is
most naturally modeled by random walks on the monoid of positive
words in ${\mathcal{LF}}_n$.

Let $d: \mathbb{N}\cup \{0\} \longrightarrow \mathbb{R}$ be a
moderating probability distribution. Following the analogy with
our constructions for a free group, we introduce  an atomic
probability measure $\tilde \mu$ on $M_n$ by assigning equal
probabilities to words of equal complexity and the total measure
$d(k)$ to the set of words of complexity $k$. This can be
interpreted as running random walks on $F_n$ of random
complexities $k$ distributed with the probability density $d(k)$.
The probability for a random walk to stop at the element $x \in
F_k$ defines a measure $\mu(x)$ on $F_k$. Obviously,
$$
\mu(x) = \sum_{\pi(w)=x} \tilde\mu(w).
$$

\subsection{How is  this new measure  related to   Kolmogorov
complexity measures?} The measure that we have just defined
should have properties very close to   Kolmogorov complexity
measures described in Section~\ref{sec:Kolmogorov}. The word $x =
\pi(w)$ is obtained from a word $w$ by cancelling   all adjacent
pairs of elements of the form $x_i X_i$, which amounts to running
a certain Turing machine on the input word $w$.  Therefore, in
view of the previous discussion, the Kolmogorov complexity $c(x)$
of $x$ differs from the Kolmogorov complexity $c(w)$ of $w$ by at
most an additive constant:
$$
c(x) \leqslant c(w) + C.
$$
Since for most words of a fixed big length $l$,   Kolmogorov
complexity is close to the length of the word, we should expect
roughly the same asymptotic behaviour
from  a random walk measure
as we do from  a measure associated with   Kolmogorov complexity.

\subsection{Comparing the length of a random walk with the geodesic
length} In  a simple case, where the complexity function on $M_n$
is just the usual length function $l(w)$, we can similarly
compare its behaviour with the behaviour of the geodesic length
$l=l_{\rm geod}(x)$ of the end point $x=\pi(w)$ of the walk on
the Cayley tree $\Gamma(F_n)$ of $F_n$, described by the word $w$.

It is easy to see that the  function $x \mapsto l_{\rm geod}(x)$
maps a random walk on $\Gamma(F_n)$ to a non-symmetric random
walk ${\mathcal{W}}_+$ on the set $\mathbb{N} \cup \{0\}$ of
nonnegative integers with reflection at $0$. In this new walk, we
make steps of length $1$; we move from the point $l=0$ to the
right with probability $1$, and, from any other point $l\ne 0$,
we move to the right with the probability $p=\frac{2n-1}{2n}$ and
to the left with the probability $q = \frac{1}{2n}$. Obviously,
the mean value of $l$ is at least the mean value of $l$ for a
random walk $\mathcal{W}$ on $\mathbb{Z}$ without reflection at
$0$; in this walk we start at $0$ and move by $1$ to the right
with the probability $p$ and by $1$ to the left with the
probability $q$. In $\mathcal{W}$, we make $m$ moves to the right
with the probability ${k \choose m}\cdot p^m\cdot q^{k-m}$ and
end up at the point $l'=m-(k-m) = 2m-k$. Thus, the random
variable $l'$ is  a linear transformation of the random variable
$m$ distributed according to the binomial distribution. Since the
expectation $E(m) = pk$, we deduce that $E(l') = 2pk-k =
\frac{n-1}{n}k$.

Therefore, for the random walk ${\mathcal{W}}_+$, the expected
value of $l$ is bounded from below by $\frac{n-1}{n}k =E(l')
\leqslant E(l)$. The upper bound $E(l) \leqslant k$ is obvious.
Thus,  the expected geodesic length $E(l_{\rm geod}(x))$ of an
element $x \in F_n$ produced by a random walk of length $k$ is
estimated as
$$
\frac{n-1}{n}k \leqslant E(l_{\rm geod}(x)) \leqslant k.
$$

\section{Behaviour of the induced measures on finite factor
groups}

In this section we try to establish possible criteria to evaluate
 how ``natural"  a given measure is, i.e.,   how    it matches 
our expectations of what the probability of hitting particular
sets should be. In general, these expectations may, of course,
vary from individual to individual, but we have some common
grounds, at least, in the case of finite sets. So, our idea is to
set up some tests to compare a given measure $\mu$ on a free group $F
= F_n $ with the induced measures on finite quotients of $F$. Most
people will probably agree that, for example, the probability for
a randomly chosen element of $F$ to have even length should be
about $1/2$. More generally, it is natural to expect that a subgroup
$H < F$  of an index $m$ would have the measure (approximately)
equal to $\frac{1}{m}$. Unfortunately, this is not the case if
$\mu(1) > 0$, indeed, if $m >> \frac{1}{\mu(1)}$ then $\mu(H) \geqslant
\mu(1) >> \frac{1}{m}$.

Moreover, most people working with probabilities on finite groups
are likely to prefer measures which are well behaved with respect
to taking ``big" finite factors $F_n/R$, for example, measures
which yield reasonable bounds for the total variance distance
\begin{equation}
\frac{1}{2}\sum_{\bar{g} \in F_n/R} \left|\mu(\bar{g}) -
\frac{1}{|F_n:R|}\right|   \label{eq:deviation}
\end{equation}
of the induced distribution  from the uniform distribution on
$F_n/R$. Again, this natural condition cannot be satisfied
basically for the same reason, namely, the value of (\ref{eq:deviation}) is
bounded from below by $\frac{1}{2}|\mu(1) - (1/m)|$, where $m = |F_n:R|$, 
and does not converge to $0$ as $m \longrightarrow \infty$.

This is a typical example of an unexcusable  dependance of a
measure on short elements. This is one of the many reasons to
believe in the following (maybe controversial) metamathematical
thesis:

\medskip
\noindent {\bf There is no particular atomic  measure on an
infinite group which would meet our expectations of sizes of
particular sets in the group}

\smallskip

A possible practical outcome of  this thesis    is either to
consider a whole family of parametric measures on a group instead
of a fixed one, or to adjust our criteria to allow a margin of an
error of approximation. The first approach was developed to some
extent in \cite{BMR}, here we make an attempt to consider the
second one.

 One may wish to exclude anomalously big
probabilities of short elements by taking  the measure of   an
element $gR$ in a finite factor group $F_n/R$ to be the
renormalized measure of the set of ``large" elements in $gR$:
\begin{equation}
\bar{\mu}_l(gR) = \frac{\mu((F_n \smallsetminus B_l) \cap gR)
}{\mu(F_n \smallsetminus B_l)},
\end{equation}
where $B_l$ is the ball of radius $l$ centered at $1$. Anyway, it
is only natural to assume that, when assessing the average case
complexity of algorithmic problems  of practical interest, we are
working with words of sufficiently large size; the bias towards
short elements can be safely ignored. A measure $\mu$ can be
accepted as ``good" if, for values of $l$ much smaller than $m$,
\begin{equation}
\frac{1}{2}\sum_{\bar{g} \in F_n/R} \left|\bar\mu_l(\bar{g}) -
\frac{1}{|F_n:R|}\right| < \frac{1}{e},
\label{eq:normalised-deviation}
\end{equation}
or is bounded by some other reasonable constant, or decreases
exponentially with  the growth of $l$:
\begin{equation}
\frac{1}{2}\sum_{\bar{g} \in F_n/R} \left|\bar\mu_l(\bar{g}) -
\frac{1}{|F_n:R|}\right| = o(e^{-cl}).
\label{eq:exponential-deviation}
\end{equation}
 The reader probably  feels already  at this point that we
are leaning  towards the classical concept of a random walk on a
group. Taking images of ``sufficiently long" elements means
allowing a sufficiently long random walk on the finite group
$F_n/R$.  It is reasonable to take the mixing time of the random
walk on the factor group $F_n/R$ with respect to the generators
$x_1R,\ldots, x_n R$ for the characteristic word length $l$ in the
criterion (\ref{eq:normalised-deviation}). Recall that the mixing
time is defined as the minimum $l_0$ such that
\begin{equation}
\|P_{l_0}-  U\| = \frac{1}{2}\sum_{\bar{g} \in F_n/R} \left|
P_{l_0}(\bar{g}) - \frac{1}{|F_n:R|}\right| < \frac{1}{e},
\label{eq:mixing-time}
\end{equation}
where $ P_l(\bar{g})$ is the probability for a random walk of
length $l$ to end up at the element $\bar{g}$, and $U$ is the
uniform distribution on $F_n/R$. From the mixing time on,  a
random walk distribution converges to the uniform distribution
exponentially fast:
\begin{equation}
\| P_{kl_0} - U\| = \frac{1}{2}\sum_{\bar{g} \in F_n/R}
\left|P_{kl_0}(\bar{g}) - \frac{1}{|F_n:R|}\right| =o(e^{-ck}).
\label{eq:exponential-mixing-time}
\end{equation}

Therefore, we can consider an alternative criterion for good
behaviour of the measure on finite factor groups:

\begin{itemize}

\item[{\bf (C1*)}] For every finite factor group $G = F_n/R$ there
exists a number $l_0$ such that
$$
\frac{1}{2}\sum_{\bar{g} \in F_n/R} \left|{\bar
\mu}_{kl_0}(\bar{g}) - \frac{1}{|F_n:R|}\right| =o(e^{-k}).
$$
\end{itemize}

Of course, the value of $l_0$ (``the mixing time'') is all
important. One should expect from a `good' measure on $F_n$ that it
is approximately the same as the mixing time of a random walk on
$F_n/R$. It is worth mentioning here, that, by a result by Pak
\cite{pak}, a random walk on  a finite group $G$ with respect to
a random generating set of size $O(\log |G|)$ mixes under $O(\log
|G|)$ steps.

It is relatively easy to see that our criterion (C1*) is met by a wide class
of measures associated with the length function on $F_n$.

\begin{theorem}
If the moderating distribution $d(k)$ satisfies a rather
natural condition   $d(k) >0$ for infinitely may values of\/ $k$
then, for every finite factor group $G = F_n/R$, the induced
measure $\bar \mu_l$ on $G$ satisfies the condition {\rm (C1*)}.
\label{th:finite-factor}
\end{theorem}

\subsection{Proof of Theorem~\ref{th:finite-factor}}
Let $|G|= m$. Every probability distribution on $G$ is a vector
$p_1,\ldots, p_m$ of non-negative real numbers subject to the
condition $p_1+\cdots+p_m = 1$. Thus, the set of all probability
distributions on $G$ is the convex hull $\Delta$ of the
distributions $E_g$ concentrated at $g \in G$, i.e.,  $E_g(h) =1 $ if
$h=g$ and $0$ otherwise. We introduce on $\Delta$ the total
variance metric
$$
\| P-Q \| = \frac{1}{2}\sum_{g\in G} |P(g) - Q(g)|.
$$
A step of a random walk induces an affine transformation
$\tau: \mathbb{R}^m \longrightarrow \mathbb{R}^m$ by the rule
$$
\tau: E_g \mapsto \frac{1}{2n}\sum_{h \mbox{ {{\scriptsize \rm
adjacent to}} }g} E_h,
$$
where the adjacency is considered in the Cayley graph of $G$.

Since all vertices of the Cayley graph of   $G$ have the  same
degree, $\tau$ fixes the uniform distribution $U = (1/m,\ldots,
1/m)$. By somewhat abusing the notation, we can move the origin
of the coordinate system in $\mathbb{R}^m$ to the point $U$. It
will be convenient to restrict the action of $\tau$ to the
subspace $\mathbb{R}^{m-1}$ spanned in $\mathbb{R}^m$ by the
polytope $\Delta$. Then the distance $\| X-U\|$ becomes a  norm
on $\mathbb{R}^{m-1}$ (which we denote by $\| X\|$;  hopefully,
this will not cause a confusion), and $\tau$ becomes a linear
operator on $\mathbb{R}^{m-1}$. Note
  that the norms $\|\tau^k(E_g)\|$  of   vertices of the
convex polytope $\tau^k(\Delta)$ are all equal. Since $\tau$ maps
the vertices $E_g$ of $\Delta$ to points of $\Delta$, it does the
same to the vertices of $\tau^k(\Delta)$. As  a result, we have
the following monotonicity property:
$$
\|\tau^{k+1}(E_g)\| \leqslant \|\tau^k(E_g)\|.
$$

 We can use again the property that
 $\tau$ maps the vertices $E_g$ of
the convex polytope $\Delta$ to prove that
  $\| \tau(E_g)\|_E < \| E_g \|_E$ for all $g\in G$, where $\|\,\|_E$
stands for the standard Euclidean norm on $\mathbb{R}^{m-1}$.
It follows that
$$
\| \tau \|_E = \max_{X\in \Delta,\, X \ne 0} \frac{\|\tau(X)\|_E}{\|
X\|_E}< 1,
$$
hence $\| \tau^k(X)\|_E \longrightarrow 0$ as $k \longrightarrow
\infty$ for every distribution $X$. Moreover, $\| \tau^k(X)\|_E$
decreases exponentially: $\| \tau^k(X)\|_E = o(c^k)$.
Since any two norms on a finite dimensional space are equivalent, we have
$$
\frac{1}{C}\|X\|_E < \|X\| < C\|X\|_E
$$
for some constant $C$, and it follows that
 $\| \tau^k(X)\| = o(c^k)$.

We have, in our new notation,
\begin{eqnarray*}
\frac{1}{2}\sum_{g\in G} \left|\bar\mu_l(g) -\frac{1}{m}\right|
= \| \bar\mu_l\|
  =  \left\| \frac{1}{\sum_{k \geqslant l} d(k)} \sum_{k \geqslant l}
 d(k) \tau^k(E_1)\right\| \\
 =  \frac{1}{\sum_{k \geqslant l} d(k)} \left\| \sum_{k \geqslant l}
 d(k) \tau^k(E_1)\right\|
  \leqslant  \frac{1}{\sum_{k \geqslant l} d(k)} \left\| \sum_{k \geqslant l}
 d(k) \tau^l(E_1)\right\| \\
  =  \|\tau^l(E_1)\|\cdot \frac{\sum_{k \geqslant l} d(k)}{\sum_{k \geqslant l}
 d(k)}  =  \|\tau^l(E_1)\|.
\end{eqnarray*}

 We see that the probabilistic distribution on
$G$ produced by  a random walk of random length $\geqslant l$
converges to the uniform distribution as fast as the probability
distribution after $l$ steps of the usual random walk does. In
particular, this proves  that the   measure $\bar \mu_l$ on $G$
satisfies the condition (C1*).

\section{The growth function and asymptotic density }

\label{sec:density}

 In this section, we discuss an approach to analyze the
 asymptotic behaviour of sets in a group via asymptotic densities
 introduced in Section \ref{se:short-inf}.

Let $G$ be a group with a complexity function $c: G \rightarrow
\mathcal{N}$ (for example, a group generated by a finite set $X$
with the  length function $l_X$).

A {\it pseudo-measure} on $G$ is a real-valued nonnegative
additive function defined on some subsets of $G$. The
pseudo-measure $\mu$ is called {\it atomic} if $\mu(S)$ is defined
for any finite subset $S$ of $G$.   Let $C_k$ and $B_k$ be, 
respectively,  the sphere and the ball of radius $k$ in $G$ with
respect to the complexity $c$.

For a set $R \subseteq F$, we define its {\it spherical asymptotic
density} with respect to $\mu$ as the following limit (if it
exists):

$$
\rho^{(s)}_{\mu}(R) = \lim_{k \rightarrow \infty} \rho^{(s)}_k(R),
$$
where
$$
\rho^{(s)}_k(R) = \frac{\mu(R \cap C_k)}{\mu(C_k)}.
$$

Similarly, we define {\it disc  asymptotic density} of $R$ as the
limit (if it exists):

$$
\rho^{(d)}_{\mu}(R) = \lim_{k \rightarrow \infty} \rho^{(d)}_k(R),
$$
where
$$
\rho^{(d)}_k(R) = \frac{\mu(R \cap B_k)}{\mu(B_k)}.
$$

One can also define the density functions above using $\lim \sup$
rather then $\lim$.

For example, if $\mu $ is the cardinality function, i.e., $\mu(A)
= |A|$,  then we obtain the standard asymptotic density functions
$\rho^{(c)}$ and $\rho^{(d)}$ on $G$.

Moreover, if $\mu$ is  $c$-invariant, i.e., \ $\mu(u) = \mu(v)$
provided $c(u) = c(v)$, then the spherical asymptotic density with
respect to $\mu$ is equal to the standard spherical asymptotic
density:
$$s\rho_{\mu} = s\rho.$$

 \begin{lemma}
 For any atomic pseudo-measure $\mu$ on $G$, the
 asymptotic densities $\rho^{(c)}$ and $\rho^{(d)}$ are  also atomic
 pseudo-measures on $G$.
\end{lemma}
 {\it Proof } is obvious.

\begin{lemma}
Let $\mu$ be a pseudo-measure on $G$. Suppose that
 $$\lim_{k \rightarrow \infty} \mu(B_k) = \infty.$$
 Then, for any subset $R \subseteq G$, if the spherical
asymptotic density $\rho^{(c)}_{\mu}(R)$ exists, then the disc
asymptotic density $\rho^{(d)}_{\mu}(R)$   exists, too,  
and
$$\rho^{(c)}_{\mu}(R) = \rho^{(d)}_{\mu}(R).$$
\end{lemma}
{\it Proof.} Let $x_k =\mu(R \cap B_k)$ and $y_k = \mu(B_k)$. Then
$y_k < y_{k+1}$ and $\lim y_k = \infty$. By Stolz's theorem
$$
\lim_{k \rightarrow \infty} \frac{x_k}{y_k} = \lim_{k \rightarrow
\infty} \frac{x_{k} - x_{k-1}}{y_{k} - y_{k-1}}.$$
 Hence
  $$ \rho^{(d)}_{\mu}(R) =
  \lim_{k\rightarrow \infty} \frac{\mu(R \cap S_k)}{\mu(S_k)} =
\rho^{(s)}_{\mu}(R),$$
as claimed.

In view of this result we will refer to the standard
(spherical or disc) densities   $\rho$.

Asymptotic densities provide a useful, though very coarse, tool to
describe behaviour of sets at infinity.
 Furthermore, there are very natural subsets which are sadly unmeasurable
 with respect to $\rho$.
To see this, consider the set $E_n$  of words of even length in a
free group  $F_n$ of rank $n$. Then $\rho(E_n)$ is easily seen to
be undefined. A way around this problem (involving generalized
summation methods for series) is discussed in \cite{BMR}.

 Moreover, the following well-known result
 shows that that  $\rho$ is not sufficiently sensitive.

\begin{theorem} {\rm (Woess \cite{woess})}
If\/ $N$ is a normal subgroup of $F_n$, $n \geqslant 2$, of
infinite index, then $\rho(N) = 0$.
\end{theorem}

Another disappointment is offered by the following result.  
 As usual, we call an element $u \in F_n$ {\it primitive} if it is part
of a free basis of   $F_n$, or, equivalently, if $\alpha(u)=x_1$
for some $\alpha \in Aut(F_n)$.

\begin{theorem} \label{2.1} Let $F_n$ be the free group of a finite rank $n
\geqslant 2$. Then:
$$\rho({\rm Pr}_n)=0,$$
where ${\rm Pr}_n$ is the set of all
primitive elements of the group $F_n$. More precisely, if $P(n,k)$ is 
the number of  primitive elements of length $k$ in $F_n, ~n \ge 3$, then 
for some constants $c_1, ~c_2$, one has   
$$c_1 \cdot (2n-3)^k \le P(n,k) \le c_2 \cdot (2n-2)^k.$$  

\end{theorem}

\smallskip

\medskip

 We see that the asymptotic density $\rho$ is not sensitive enough
in measuring sets in $F_n$.
It would be very interesting therefore to check whether
probabilistic results of, say,
 \cite{Arzh}, \cite{Fr}, \cite{Olsh}, that are based on the
asymptotic density  $\rho$, will hold upon replacing $\rho$ with a
more adequate  measuring tool.

\subsection{Proof of Theorem~\ref{2.1}} Our proof  is based on the fact that
the Whitehead graph of any primitive element of length $>2$ has
either an isolated edge or a cut vertex, i.e., a vertex  that,
having been removed from the graph
 together with all  incident edges, increases the number of
connected components of the graph. Recall that the Whitehead
graph $Wh(u)$  of a word $~u$ is obtained as follows. The
vertices of this graph correspond  to the elements of the free
generating set $X$ and their inverses. If the  word $~u$ has a
subword  $x_i x_j$, then there is an edge in $Wh(u)$ that
connects the vertex $x_i$ to the vertex $x_j^{-1}$; ~if $~u$ has
a subword  $x_i x_j^{-1}$, then there is an edge that connects
$x_i$ to $x_j$, etc. ~We note that usually, there is one more
edge (the external edge) included in the definition of the
Whitehead  graph: this is the edge that connects the vertex
corresponding to the last letter of $~u$, to the vertex
corresponding to the inverse of the first letter.

 Assume first that the Whitehead graph of $~u$ has a cut vertex.
We are going to show  that the number of elements $u$ of length
$k$ with this property  is no bigger than $C \cdot (2n-2)^{k-1}$,
where $C=C(n)$ is a constant.

 Let  $Wh(u)$ be a disjoint union of two graphs, $\Gamma_1$
 and $\Gamma_2$, complemented by  a (cut)  vertex $A$ together
with all incident edges. Let $n_1 \geqslant 1$  and $n_2 \geqslant
1$ be the number of vertices in $\Gamma_1$  and $\Gamma_2$,
respectively. Then, in particular, $n_1+n_2=2n-1$. Let
$m=\min(n_1, n_2)$.

The first letter of $u$ can be any of the $2n$ possible ones.
 For the following letter however we  have no more than $2n-m$
 possibilities since, for example, if the first letter, call it
$x_i$,  corresponds to a vertex from $\Gamma_1$, then the
following letter,
   call it $x_j$, cannot be such
  that the  vertex corresponding to $x_j^{-1}$ belongs to $\Gamma_2$,
because otherwise, there would be an edge in $Wh(u)$  that
connects $\Gamma_1$ to $\Gamma_2$ directly, not through $A$,
i.e., $A$ would not be a cut vertex.

 Thus, there are only $2n-n_2 \leqslant 2n-m$  possibilities for the
following letter in $u$. The same argument applies to every letter
 in $u$, starting with the second one; therefore, the total number
of possibilities (corresponding to a particular choice of $n_1$
and $n_2$) is no bigger than $C \cdot (2n-m)^{k-1}$, where
$C=C(n)$ is a constant.

 Now consider two cases:
\begin{itemize}
\item[{\bf (i)}] $m \geqslant2$.
\end{itemize}
   In this case, the number in question is
bounded by   $C \cdot (2n-2)^{k-1}$.
\begin{itemize}
\item[{\bf (ii)}] $m =1$.
\end{itemize}
  In this case, one of the graphs, say, $\Gamma_2$, consists of a single
vertex (call it $x_1$ for notational convenience) connected only
to the cut vertex (call it $x_2$); in particular, the  vertex
$x_1$ is  a terminal vertex of the graph $Wh(u)$, and, whenever
the letter $x_1$ occurs in $u$, it is followed by $x_2^{-1}$. Let
$q \geqslant1$ be the number of occurrences of $x_1$ in $u$.

 Apply the automorphism $\phi: x_1 \mapsto x_1 x_2, ~x_i \mapsto x_i, ~i \geqslant2$,
~to the element $u$. Then $|\phi(u)| = |u|-q < |u|$. Now
$\phi(u)$ is a primitive element, too; hence the whole argument
above is applicable to $\phi(u)$. In particular, if, in the
notation above, $m \geqslant 2$ for this element $\phi(u)$, then,
as we have just proved, the number of those elements  is
 bounded by $C \cdot (2n-2)^{k-1-q}$ for some $C=C(n)$. This number is also
equal to the number of elements $u$ of the type we are
considering now (for a particular  $q$) because of the one-to-one
correspondence between elements $u$ and $\phi(u)$.

 If $m =1$ for the element $\phi(u)$, then we use the same trick again, until
we get to a primitive element with $m \geqslant 2$. In any case,
the number of primitive elements of length $k$, with $m =1$,  and
with $q \geqslant1$
 occurrences of $x_1$, is  bounded by $C \cdot (2n-2)^{k-1-q}$ for some $C=C(n)$.

 Now we have to sum  up for all possible values of $q$. Note that $q$ cannot be
equal to $k$ since $x_1^k$ is not a primitive element; also, $q$
cannot be equal to $k-1$ since   in the Whitehead  graph of
$x_1^{k-1}x_2$, the vertex corresponding to $x_1$ is not a
terminal vertex. Hence, we have

$$\sum_{q=1}^{k-2} C \cdot (2n-2)^{k-1-q} = C \cdot \frac{(2n-2)^{k-1} - 1}{2n-3}
 \leqslant  C \cdot (2n-2)^{k-1}.$$

 Finally, we have to multiply this number by the number of ways we can choose
two vertices (terminal and cut) out of $2n$, i.e., by $(2n-1)n$,
but this does not     change the type of the bound. Therefore,
 we get here the same bound as we got in the case $m \geqslant 2$.

\smallskip

  Thus,  summing up for all possible values of $n_1 \geqslant 1$ and  $n_2 \geqslant 1$
  such
  that $n_1+n_2=2n-1$, we see that  the total number of possibilities in (i)
 and (ii) is bounded by $C \cdot (2n-2)^k$ for some $C=C(n)$.
\smallskip

 Finally,   we address the remaining case, where the Whitehead
graph of $u$ has  an isolated edge. In that case, some cyclic
permutation of $u$ must be of the form $x_i^{\pm 1} x_j^{\pm 1}
u_1$, where $u_1$ does not depend on  $x_i$,   $x_j$, and $j$
does not have to be different from $i$. The number of elements
with this property is easily seen to be bounded by $C \cdot k
\cdot (2n-3)^{k-1}$ ~for some constant $C=C(n)$.
\smallskip

Therefore, the ratio of the number of primitive elements
 of length $k$ to the number of all elements
 of length $k$ is no more than $C \cdot \frac{(2n-2)^k}{2n (2n-1)^{k-1}} =
 C' \cdot
 (\frac{2n-2}{2n-1})^k$, ~where $C'=C'(n)$ is a constant.
This ratio obviously tends to 0, and, moreover, well-known
properties of a geometric series now imply that the ratio of the
number of primitive elements
 of length $\leqslant k$ to the number of all elements
 of length $\leqslant k$ tends to 0, too.

 Finally, we note that the lower bound $c_1 \cdot (2n-3)^k \le P(n,k)$ 
is obvious because every element of the form $x_1 \cdot u(x_2, ..., x_{n-1})$
is primitive. 

\medskip

  Just to complete the picture, we also mention here the
 bounds for the number of primitive elements of length $k$ in  $F_2$:

 \begin{proposition} \label{2.2} {\rm (\cite{MSh})} The number of primitive
elements of length\/ $k$ in the group $F_2$ is:

\begin{itemize}
\item[{\bf (a)}] more than\/ $\frac{4}{\sqrt 3} \cdot (\sqrt 3)^k$
if\/ $k$ is odd.

\item[{\bf (b)}] more than\/ $\frac{4}{3} \cdot (\sqrt 3)^k$ if\/ $k$ is
even.

\end{itemize}
\end{proposition}

Informally speaking, ``most" primitive elements  in $F_2$ are
conjugates of primitive elements  of smaller length. This is not
the case in $F_n$ for $n > 2$, where ``most"  primitive elements
are of the form $u \cdot  x_i^{\pm 1} \cdot v$, where  $u, v$
are   arbitrary elements that do not depend on $x_i$.
\medskip

\subsection{The rate of convergence of the asymptotic density}
\label{sec:density-2}

Let $S$ be a subset of a free group $F_n$ of rank $n \geqslant 2$.
As we have mentioned in the Introduction, the growth  rate of the
set $S$ is  defined as the limit
$$\gamma(S) = \limsup_{k \to \infty}\sqrt[k] {\rho_k(S)},$$
where
$$\rho_k(S) = \frac{|S\cap  B_k|}{|B_k|}$$
are the disc relative frequencies of $S$. The growth rate
$\gamma(S)$ is  gauging the speed of convergence to zero  of the
sequence $\rho_k(S)$. Indeed,  the  standard results from calculus
show that
 if $\gamma(S) < 1$, then $\rho(S) = 0$, and, moreover,
the sequence $\rho_k(S)$ converges to $0$ exponentially fast.

A classical theorem by Grigorchuk \cite{grigorchuk} states  that
if $N \le F_n$ is a normal subgroup 
 and $F_n/N$ is not an amenable group, then
$\gamma(N) < 1$. This gives us many easy examples of normal
subgroups $N$ with the exponential speed of convergence of the
asymptotic density. For example, this happens every time when the
factor group $F_n/N$ contains a non-abelian free group.

It would be interesting to have a look at the other end of the
spectrum and estimate the speed of convergence of the sequence
$\rho_k(N)$ for an obviously ``big" subgroup $N \le F_n$. The
following result is one of the very few instances where we have
concrete information:

\begin{theorem} {\rm (Sharp \cite{sharp})} If $n \geqslant 2$, then the spherical
relative frequencies
$$
\rho^{(s)}_k = \frac{|[F_n,F_n] \cap S_k|}{|S_k|}
$$
of words of length $k$ from the derived subgroup $[F_n,F_n]$ of
the free group $F_n$ asymptotically behave as
$$
\rho^{(s)}_k \sim \left\{\begin{array}{cl}
0 & \quad \hbox{ if} \quad k \;\hbox{ is odd},\\
 \frac{C}{k^{n/2}} & \quad \hbox{ if} \quad k \;\hbox{ is even.}
 \end{array}\right.
$$
\end{theorem}

\end{document}